\documentclass[11pt, reqno, english]{amsart}

\usepackage{amscd,amsmath}
\usepackage{mathrsfs}
\usepackage{amsfonts}
\usepackage{amssymb}
\usepackage{enumerate}
\usepackage{setspace}
\usepackage{color}

\usepackage[colorlinks=true,
            linkcolor=red,
            urlcolor=blue,
            citecolor=gray]{hyperref}

\usepackage{amsthm}
\usepackage{amssymb,latexsym,graphics,enumerate}
\usepackage[mathscr]{eucal}
\usepackage{amsmath,amsfonts,amsthm,amssymb}
\usepackage[left=1in,right=1in,top=1in,bottom=1in]{geometry}
\usepackage{tikz} 

\newcommand{\rr}{\mathbb{R}}

\newcommand{\be}{\begin{eqnarray*}}
\newcommand{\bel}{\begin{eqnarray}}
\newcommand{\ee}{\end{eqnarray*}}
\newcommand{\eel}{\end{eqnarray}}
\newcommand{\ba}{\begin{aligned}}
\newcommand{\ea}{\end{aligned}}

\newcommand{\al}{\alpha}

\newcommand{\na}{\nabla}

\newcommand{\pa}{\partial}
\newcommand{\bs}{\backslash}

\newtheorem{thm}{Theorem}[section]

\newtheorem{rmk}{Remark}[section]

\renewcommand{\geq}{\geqslant}

\renewcommand{\leq}{\leqslant}

\def\d{{\sf{d}}}
\def\e{{\sf{e}}}

\def\h{{\sf{h}}}

\DeclareMathOperator{\trace}{Tr}
\DeclareMathOperator{\diver}{div}

\def\lam{\lambda}
\def\hf{\frac{1}{2}}

\def\dy{\text{d}y}

\def\etaS{\eta_{{}_S}}
\def\etaSz{\eta_{{}_{S_0}}}
\def\etaM{\eta_{{}_M}}
\def\cm{c_{\text{min}}}

\def\u{{\mathbf u}}
\def\v{{\mathbf v}}
\def\w{{\mathbf w}}
\def\bs{{\mathbf s}}
\def\br{{\mathbf r}}
\def\bell{{\boldsymbol \ell}}
\def\diam{D}
\def\deg{\text{deg}}
\def\supp{\text{supp}}
\def\suppr{\supp\{\rho(\cdot,t)\}}

\numberwithin{equation}{section}

\begin{document}

\title[Global Regularity of Two-Dimensional  Flocking Hydrodynamics]{Global regularity of\\two-dimensional  flocking hydrodynamics}


\author{Siming He}
\address{Department of Mathematics and Center for Scientific Computation and Mathematical Modeling (CSCAMM), University of Maryland, College Park}
\email{simhe@cscamm.umd.edu}

\author{Eitan Tadmor}
\address{Department of Mathematics, Center for Scientific Computation and Mathematical Modeling (CSCAMM), and Institute for Physical Sciences \& Technology (IPST), University of Maryland, College Park\\
Current address: ETH Institute for Theoretical Studies,
ETH-Z\"urich, 8092 Z\"urich, Switzerland}
\email{tadmor@cscamm.umd.edu}

\date{\today}

\subjclass{92D25, 35Q35, 76N10}

\keywords{flocking, alignment, hydrodynamics, regularity, critical thresholds.}

\thanks{\textbf{Acknowledgment.} Research was supported in part by NSF grants DMS16-13911, RNMS11-07444 (KI-Net) and ONR grant N00014-1512094. We thank  the ETH Institute for Theoretical Studies (ETH-ITS) for the support and hospitality.}
\date{\today}

\begin{abstract}
We study the systems of Euler equations which 
  arise from agent-based dynamics driven by velocity \emph{alignment}. 
  It is known that smooth solutions of such systems must flock, namely --- the large time behavior of the velocity field  approaches a limiting ``flocking'' velocity. To address the question of global regularity, we derive sharp critical thresholds in the phase space of initial configuration which characterize the global regularity and hence flocking behavior of such \emph{two-dimensional} systems. Specifically, we prove for that a large class of \emph{sub-critical}  initial conditions  such that the initial divergence is ``not too negative'' and the initial spectral gap is ``not too large'',  global regularity persists for all time.  
\end{abstract}

\ifx
\begin{altabstract}
Nous \'{e}tudions les syst\'{e}mes des \'{e}quations d'Euler qui
  r\'{e}sultent de dynamiques d'\emph{alignement} entre agents.
Il a \'et\'e prouv\'e  que pour des solutions r\'eguli\`eres de tels syst\'{e}mes, en temps grand,  le champ de vitesse s'approche d'une vitesse limite uniforme. Nous identifions des seuils critiques dans l'espace de phase de la configuration initiale qui caract\'{e}risent la régularit\'{e} globale et donc le comportement en temps grand de tels syst\`{e}mes \emph{bidimensionnels}. 
Plus pr\'{e}cis\'{e}ment, nous prouvons que pour une classe assez large de conditions initiales \emph{sous-critiques} telles que la divergence initial n'est ``pas trop n\'{e}gative'' et l'\'{e}cart spectral initial n'est ``pas trop grand'', la r\'{e}gularit\'{e} globale reste vraie en temps grand.
\end{altabstract}
\fi

\maketitle
\setcounter{tocdepth}{1}
\tableofcontents
\section{Flocking hydrodynamics}
We consider the system of Eulerian dynamics where the density $\rho(x,t)$ and velocity field $\u(x,t) = (u_1, \ldots u_n): {\mathbb R}^n \times {\mathbb R}_+ \mapsto {\mathbb R}^n$ are driven by nonlocal alignment forcing,

\bel\label{The basic equation}
\left\{
\ba
&\rho_t+\na\cdot(\rho \u)=0,\\
&\u_t+\u\cdot\na \u=\int a(x,y,t)(\u(y,t)-\u(x,t))\rho(y,t)\dy
\ea \right\} \quad  (x,t) \in \rr^n\times \rr_+.
\eel

A solution $(\rho,\u)$ is sought subject to the compactly supported initial density $\rho(x,0)=\rho_0(x)\in L_+^1(\rr^n)$ and uniformly bounded initial velocity $\u(x,0)=\u_0(x)\in W^{1,\infty}(\rr^n)$. The alignment forcing on the right hand side of (\ref{The basic equation}) involves the non-negative  interaction kernel $a(x,y,t)$. \newline
Such systems arise as 
macroscopic realization of agent-based dynamics which describes the collective motion of $N$ agents, each of which adjusts its velocity to a \emph{weighted average} of  velocities of its neighbors 
\bel\label{agent based model}
\left\{\begin{array}{rr}\ba&\dot{x}_i=\v_i\\
&\dot{\v}_i=\frac{1}{\deg_i}\sum_{j=1}^N \phi(|x_i-x_j|)(\v_j-\v_i)\ea\end{array}\right.
\eel
Here, the weighted average of the right of \eqref{agent based model} is dictated by  influence function $\phi(\cdot)$ which is assumed to be decreasing, and $\deg_i$ is a weighting normalization factor.
Different  agent based models employ different $\deg_i$'s, e.g., \cite{CCP2017}. We focus here on two such models. The Cucker-Smale (CS) model \cite{CS2007} sets a uniform averaging $\deg_i\equiv N$ which leads to the \emph{symmetric} interaction kernel $a(x,y)=\phi(|x-y|)$. The Motsch-Tadmor (MT) model \cite{MT2011} uses an \emph{adaptive} normalization  $\deg_i=
\sum_j \phi(|x_i-x_j|)$ which leads to $\displaystyle a(x,y,t)=\frac{\phi(|x-y|)}{(\phi*\rho)(x,t)}$. The kernel is non-symmetric but normalized such that  $\int a(x,y,t)\rho(y,t)\dy=1$.
The dynamics of \eqref{agent based model} can be described in terms of the empirical distribution
$f(x,\v,t):=\frac{1}{N}\sum_j \delta_{x-x_j(t)}\otimes \delta_{\v-\v_j(t)}$.
For large crowds of $N$ agents, $N\gg1$,  a limiting distribution of the approximate  form  $f(x,\v,t) \approx \rho(x,t)\delta(\v-\u(x,t))$   is captured by the first two  velocity moments, namely -- the density $\rho:= \langle f(x,\v,t)\rangle$ and momentum $\rho \u :=\langle \v f(x,\v,t)\rangle$ satisfy the \emph{conservative} system \cite{HT2008,CCR2009,CFRT2010,MOA2010}
\begin{equation}\label{eq:CS+MT}
\hspace*{-0.3cm}\left\{
\ba
\rho_t+\na\cdot(\rho \u)&=0\\
(\rho\u)_t+\na (\rho \u\otimes \u)&=\frac{\al(x,t)}{(\phi*\rho)(x,t)}\int \phi(|x-y|)(\u(y,t)-\u(x,t))\rho(x,t)\rho(y,t)\dy.
\ea \right.
\end{equation}
Here $\al(x,t)$ is the amplitude of alignment, $\al(x,t)=(\phi*\rho)(x,t)$ in the case of CS model, and $\al(x,t)\equiv 1$ in MT model.
When classical solutions of these equations are restricted to the support of $\rho(\cdot,t)$, one ends with  the equivalent system \eqref{The basic equation} with $a(x,y,t)=\al(x,t)\phi(|x-y|)/(\phi*\rho)(x,t)$, namely
\begin{equation}\label{eq:non-vacuum}
\left\{
\ba
\rho_t+\na\cdot(\rho \u)&=0,\\
\u_t+ \u\cdot\na \u&=\frac{\al(x,t)}{(\phi*\rho)(x,t)}\int \phi(|x-y|)(\u(y,t)-\u(x,t))\rho(y,t)\dy.
\ea \right.
\end{equation}
 Since the alignment forcing on the right is non-local, dictated by the support of $\phi$, it acts even   within the vacuum region where  $ \text{dist}\{x,\suppr\}>0$, and  \eqref{eq:non-vacuum} extends \emph{throughout} $\rr^n$. We elaborate on this issue in \S\ref{sec:vacuum} below.
 
We note that the dynamics of both models can be interpreted in terms of  the mean velocity $\overline{\u}(x,t)$
\be\label{eq:fhd}
\u_t+\u\cdot\na \u=\al(x,t)\big(\overline{\u}(x,t)-\u(x,t)\big), \qquad \overline{\u}(x,t):=\frac{\phi*(\rho \u)(x,t)}{(\phi*\rho)(x,t)}.
\ee
 This formulation reveals that  system  \eqref{eq:non-vacuum} (and in its general form \eqref{The basic equation}) is dynamically aligned towards the mean $\overline{u}(x,t)$, and its large time behavior is expected to approach  a constant limiting velocity. This is the flocking hydrodynamics alluded to in the title, where a finite-size of non-vacuum  state is approaching a limiting velocity as $t \rightarrow \infty$.  Specifically, the dynamics can be characterized in terms of the diameters
\be
\diam(t):=\sup_{x,y\in \suppr}|x-y|, \qquad \
V(t):=\sup_{x,y\in \suppr}|\u(x,t)-\u(y,t)|.
\ee
The system \eqref{The basic equation} converges to a flock  if there exists a finite $D$ such that
\bel\label{S, V control}
\sup_{t\geq 0}\diam(t)\leq \diam_\infty \ \ \text{ and } \ \ V(t)\stackrel{t\rightarrow \infty}{\longrightarrow} 0.
\eel
This corresponds to the flocking behavior at the level of agent-based description \cite{HT2008}, \cite[definition 1.1]{MT2011} where a cohesive flock of a finite diameter $\max_{i,j}|x_i(t)-x_j(t)|\leq \diam_\infty<\infty$, is approaching a limiting velocity, $\max_{i,j}|\v_i(t)-\v_j(t)|\rightarrow 0$ as $t\rightarrow \infty$.

\subsection{Strong solutions must flock} In this work we focus on the case where $\phi$ is \emph{global}. Since the agent based model (\ref{agent based model}) exhibit flocking behavior in this case, \cite{MT2014}, it is natural to to expect a similar result for its macroscopic description (\ref{eq:non-vacuum}). This is the content of the following theorem.

\begin{thm}[Strong solutions must flock \cite{TT2014}]\label{Strong solutions must flock} Let $(\rho(\cdot,t),\u(\cdot,t)) \in (L^\infty\cap L^1) \times W^{1,\infty}$ be a global strong solution of the system (\ref{eq:non-vacuum}) subject to a compactly supported initial density $\rho_0=\rho(\cdot,0)\geq 0$ and bounded initial velocity $\u_0=\u(\cdot,0)\in W^{1,\infty}$. Assume that a monotonically decreasing influence function $\phi \leq \phi(0)=1$ is global in the sense that\footnote{We let $|\cdot|_p$ denote the usual $L^p$ norm.}
\bel\label{phi condition}
V_0 < m_0 \int_{\diam_0}^{\infty} \phi(r)dr, \qquad  m_0:=|\rho_0|_1,
\eel
where  $\diam_0$ and $V_0$ are the initial diameters of non-vacuum density and velocity. 
Then $(\rho, \u)$ converges to a flock at exponential rate, namely --- the support of $\rho(\cdot,t)$ remains within a finite diameter $\diam_\infty$ whose existence follows from assumption \eqref{phi condition}
\begin{subequations}\label{eqs:global}
\begin{equation}\label{eq:Dinf}
 \sup_{t\geq 0}\diam(t) \leq \diam_\infty \quad \text{where} \quad m_0\int_{\diam_0}^{\diam_\infty} \phi(s)\text{d}s = V_0,
 \end{equation}
  and  
\bel\label{S, V exp-control}
V(t) \leq V_0 e^{-\kappa t} \longrightarrow 0, \qquad  \kappa:=\left\{\begin{array}{ll}m_0\phi_\infty, & \mathrm{CS} \ \ \text{model},\\ \phi_\infty, & \mathrm{MT} \  \text{model},\end{array}\right. \ \ \phi_\infty:=\phi(\diam_\infty).
\eel
\end{subequations}
In  particular, if $|\phi|_1=\infty$ then there is an \emph{unconditional flocking}
in the sense that  \eqref{eqs:global} holds for all finite $V_0$.
\end{thm}
For the sake of completeness we provide below an alternative derivation of the exponential alignment in  \eqref{eqs:global}, as an a priori bound instead of the ``propagation along characteristics'' argument in \cite[Theorem 2.1]{TT2014}. To this end, we extend the scalar argument in \cite[Lemma 1.1]{ST2017} to general systems using a projection argument employed in \cite[Theorem 2.3]{MT2014}. Fix an arbitrary $\w\in \mathbb{R}^n$ and project the CS model \eqref{eq:non-vacuum} on $\w$ to find
\[
(\partial_t+\u\cdot\nabla)\langle\u(x,t),\w\rangle =\int \phi(|x-y|)\Big(\langle \u(y,t),\w\rangle -\langle\u(x,t),\w\rangle\Big)\rho(y,t)\text{d}y.
\] 
It follows that $\displaystyle u_+(t):=\max_{x\in \suppr} \langle \u(x,t),\w\rangle$ satisfies
\[
\begin{split}
\frac{\text{d}}{\text{d}t}u_+ & =\int \phi(|x_+-y|)\Big(\langle \u(y,t),\w\rangle -\langle\u(x_+,t),\w\rangle\Big)\rho(y,t)\text{d}y \\
& \leq 
\min_{x,y\in \text{supp}\, \rho(\cdot,t)} \phi(|x-y|) \int \Big(\langle \u(y,t),\w\rangle -\langle\u(x_+,t),\w\rangle\Big)\rho(y,t)\text{d}y
\end{split}
\]
Similarly, we have the lower bound on $\displaystyle  u_-(t):=\mathop{\min}_{x\in \suppr} \langle \u(x,t),\w\rangle$
\[
\frac{\text{d}}{\text{d}t}u_-  \geq 
\min_{x,y\in \suppr} \phi(|x-y|) \int \Big(\langle \u(y,t),\w\rangle -\langle\u(x_-,t),\w\rangle\Big)\rho(y,t)\text{d}y
\]
The difference of the last two inequalities  implies 
\[
\frac{\text{d}}{\text{d}t}|u_+(t)-u_-(t)|\leq -\phi(D_\infty) m_0 |u_+(t)-u_-(t)|, \qquad \phi(D_\infty)=\min_{x,y\in \suppr} \phi(|x-y|).
\]
It follows that the CS velocity diameter, $\displaystyle V(t)=\sup_{|\w|=1} |u_+(t)-u_-(t)|$, satisfies the bound  \eqref{S, V exp-control}  with $\kappa=m_0\phi_\infty$. The same argument  follows for  MT model with $\kappa = \phi_\infty$, independently of $m_0$.

\subsection{Critical thresholds} Theorem \ref{Strong solutions must flock} raises the problem whether solutions of the hydrodynamic model (\ref{eq:non-vacuum}) remain smooth for all time.  This question was  addressed in \cite{TT2014, CCTT2016}, proving that the compactly supported initial data stay below certain critical threshold in configuration space  then initial  smoothness propagates and as a result, the corresponding strong solutions will flock. 
Recall the finite-time blow-up  of compactly supported density in the presence of  \emph{local} pressure \cite{Si1985,LY1997} and even in the presence of global Poisson forcing \cite{Ma1992}.
In both cases, a positive lower-bound on the (potential of) the forcing --- the pressure, Poisson, etc, over the \emph{finite} $\suppr$ leads to finite time blow up.
In contrast, here the non-local character of the influence function $\phi$ guarantees global regularity, at least for sub-critical initial data. This type of conditional regularity for Eulerian dynamics depending on a \emph{critical threshold} in configuration space,  was advocated in a series of papers \cite{ELT2001,LT2002,LT2003,LT2004,HT2008,LL2013}. Here, we pursue this approach to derive sharp critical thresholds  for propagation of regularity of the \emph{two-dimensional} flocking hydrodynamics.

\subsection{Vacuum and the finite horizon alignment}\label{sec:vacuum} According to \eqref{phi condition}, if the influence function is global in the sense that $\displaystyle \int^\infty \!\!\phi(r)dr=\infty$, then the alignment dynamics \eqref{eq:non-vacuum} admits \emph{unconditional} flocking in the sense that \eqref{eqs:global} holds for all $V_0$'s. This holds for both the symmetric CS model and non-symmetric MT model \cite[proposition 2.9]{MT2014}. In this case, 
 alignment in \eqref{eq:non-vacuum} is active \emph{throughout} $\rr^n$, inside and outside $\suppr$. Indeed, one has a global lower-bound on the action of alignment for all $x\in \rr^n$, \cite[proposition 6.1]{TT2014}
\[
(\phi*\rho)(x,t) \geq m_0 \phi(d(x,t)+D_\infty) >0, \qquad d(x,t)=\text{dist}\{x, \suppr\}
\]
The flocking behavior of such a  global approach was pursued in \cite{TT2014}.\newline
Another possible approach to study \eqref{eq:non-vacuum} is to focus on a specific initial configuration with finite velocity variation $V_0<\infty$.
Then, since  $\suppr$ cannot grow beyond a maximal  diameter of size $D_\infty$ dictated by \eqref{eq:Dinf}, it follows that the alignment term on the right of the underlying conservative formulation \eqref{eq:CS+MT}, 
\[
\phi(|x-y|)(\u(y,t)-\u(x,t)) \rho(x,t)\rho(y,t)\equiv 0, \qquad |x-y|> D_\infty,
\]
independently of  the values of $\{\phi(r),\ r>D_\infty\}$. Alternatively, we can fix a compactly support influence function $\phi$ and view \eqref{eq:Dinf} as a restriction on initial velocities whose variation is ``not too large'', so that they lead to flocking. With either one of these two points of view,  the values of $\phi(r)$ for $r>D_\infty$ play no role in the dynamics.
We therefore may set $\phi(r)_{|r>D_\infty} \equiv 0$ which in turn sets a \emph{finite horizon} on the action of  alignment. Namely, the alignment in \eqref{eq:non-vacuum} is still active in the vacuous annulus \emph{outside} $\suppr$,
\[
A(t):=\{x \ | \ 0< \text{dist} \{x, \suppr\} <D_\infty\},
\]
and \eqref{eq:non-vacuum} applies in $\suppr\cup A(t)$,
\begin{subequations}\label{eqs:horizon}
\begin{equation}\label{eq:final}
\hspace*{-0.6cm}\left\{
\ba
\rho_t+\na\cdot(\rho \u)&=0,\\
\u_t+ \u\cdot\na \u&=\frac{\al(x,t)}{\phi*\rho}\int \phi(|x-y|)(\u(y)-\u(x))\rho(y)\dy   
 \ea \right\} \  \text{dist} \{x, \suppr\} <D_\infty.
\end{equation}
However, since $\phi(|x-y|)\rho(y)$ is supported for $y$'s in the intersection 
$y\in Y_x(t):= \suppr\cap B_{D_\infty}(x)$, it implies the alignment bound
\[
\left|\int \phi(|x-y|)(\u(y,t)-\u(x,t)) \rho(y,t)dy\right| \leq  V(t) \cdot |\rho(\cdot,t)|_\infty\times\int_{y\in Y_x(t)}\phi(|x-y|)dy.
\]
It follows that the alignment on the right of \eqref{eq:final} approaches zero, as $x \in A(t)$ approaches the ``horizon'' 
boundary $\text{dist} \{x, \suppr\} =D_\infty$ and $\text{vol}(Y_x(t)) \rightarrow 0$. In particular, 
$(\phi*\rho)(x,t) \equiv 0$ beyond the horizon $\text{dist} \{x, \suppr\} >D_\infty$, where the momentum equation is reduced to inviscid pressureless equations,
$\u_t+\u\cdot\nabla\u=0$.
Accordingly, \eqref{eq:final} can be complemented with  \emph{constant} far-field boundary conditions, in agreement with \cite[Remarks 2.8 \& 6.6]{TT2014}, 
\begin{equation}\label{eq:BC}
\u(x,t)\equiv \u_\infty, \quad \text{for} \ \ \text{dist} \{x, \suppr\} >D_\infty.
\end{equation}
\end{subequations}

\section{Cucker-Smale hydrodynamics}
\subsection{Global regularity}
We begin by recalling the one-dimensional Cucker-Smale model for $(\rho,u): (\rr,\rr_+) \mapsto (\rr_+,\rr)$,
\bel\label{CS1D}
\left\{\ba
&\rho_t+(\rho u)_x=0,\\
&u_t+uu_x=\int_{\rr} \phi(|x-y|)(u(y,t)-u(x,t))\rho(y,t)\dy\\
\ea\right.\qquad (x,t) \in (\rr,\rr_+).
\eel
In \cite{CCTT2016} it was proved that  (\ref{CS1D}) has a global classical solution if and only if the initial data satisfies 
\begin{equation}\label{eq:CT1D}
\pa_x u_0(x)\geq -(\phi*\rho_0)(x), \ \mbox{ for all} \ x \ \ \in \rr.
\end{equation}
Condition \eqref{eq:CT1D} separates the space of initial configurations into two distinct  regimes: a sub-critical regime of initial data satisfying  $\pa_x u_0(x)\geq -\phi*\rho_0(x),\forall x\in supp(\rho_0)$, which guarantee  \emph{global} smooth solutions; and a supercritical regime of initial conditions  such that $\pa_x u_0(x_0)\leq -\phi*\rho_0(x_0)$ for some $x_0\in \rr$, which leads to a finite time blowup. 
This is a typical one-dimensional example for the {critical threshold} behavior. Condition \eqref{eq:CT1D} provides a sharp improvements to the earlier critical threshold results in \cite{ST1992,LT2001,TT2014}. Recent results in \cite{ST2016,DKRT2017} prove the global regularity of \eqref{CS1D} for singular kernels $\phi(|x|)=|x|^{-(1+\alpha)}$ for $\alpha\in (0,2)$ \emph{independent} of any finite critical threshold. Singularity helps!. 

A first attempt to extend the study of critical threshold to the \emph{two-dimensional} CS model  was derived in \cite{TT2014}. Here, we improve this result with a simplified derivation  of a sharper critical threshold condition, leading to alignment decay of order $e^{-\kappa t}$. We recall  \eqref{S, V exp-control} which set $\kappa=m_0\phi_\infty$ in the present case of CS model. 

\begin{thm}[Critical threshold for 2D Cucker-Smale hydrodynamics]\label{thm:CS}
\noindent
Consider the two-dimensional CS model 
\bel\label{CS}
\left\{\ba
&\rho_t+\na\cdot(\rho \u)=0,\\
&\u_t+\u\cdot\na \u=\int \phi(|x-y|)(\u(y,t)-\u(x,t))\rho(y,t)\dy\\
\ea\right\} \quad x\in \rr^2, t\geq 0,
\eel
subject to initial conditions, $(\rho_0, \u_0)\in (L_+^1(\rr^2), W^{1,\infty}(\rr^2))$, 
with compactly supported density, $D_0<\infty$, and such that the variation of the initial velocity satisfies the strengthened bound 
\bel\label{V_0}
V_0\leq m_0\cdot\min\left\{|\phi|_{1},\frac{\phi^2_\infty}{4|\phi'|_\infty}\right\} , \qquad
V_0=\mathop{\max}_{x,y\in \text{supp}(\rho_0)}|\u_0(x)-\u_0(y))|, \quad \phi_\infty=\phi(\diam_\infty).
\eel

\noindent
Assume that the following critical threshold condition holds.\newline
 (i) The initial velocity divergence satisfies
\bel\label{eq:CT2D}
\diver{\u_0}(x) \geq -\phi*\rho_0(x) \quad \mbox{for all} \quad x\in \rr^2.
\eel
(ii) Let $S=\hf\{(\pa_ju_i+\pa_i u_j)\}$ denote the symmetric part of the velocity gradient with  eigenvalues  $\mu_i=\mu_i(S)$. Then the initial \emph{spectral gap} $\etaSz:=\mu_2(S_0)-\mu_1(S_0)$ is bounded
\begin{equation}\label{eq:etaCT}
\max_x \big|\etaSz(x)\big| \leq \frac{1}{2}m_0 \phi_\infty, \qquad \etaS=\mu_2(S(x,t))-\mu_1(S(x,t)).
\end{equation}  
Then the class of such sub-critical initial conditions \eqref{eq:CT2D},\eqref{eq:etaCT} admit a classical solution\newline $(\rho(\cdot,t),\u(\cdot,t))\in C(\rr^+;L^\infty\cap L^1(\rr^2 ))\times C(\rr^+;\dot{W}^{1,\infty}(\rr^2))$ with large time hydrodynamics flocking behavior \eqref{S, V exp-control}, $\displaystyle \mathop{\max}_{x,y\in \text{supp}(\rho(\cdot,t))} |\u(x,t)-\u(y,t)| \lesssim e^{-\kappa t}$.
 \end{thm}

\noindent
Before turning to the proof of theorem \ref{thm:CS}, we comment on its assumptions.
\begin{rmk}[on the critical threshold \eqref{eq:CT2D},\eqref{eq:etaCT}] Theorem \ref{thm:CS} recovers the one-dimensional critical threshold \eqref{eq:CT1D}. It amplifies the same theme of  critical threshold required for global regularity of other \emph{two-dimensional} Eulerian dynamics found in restricted Euler-Poisson \cite{LT2003}, rotational Euler \cite{LT2004},..., namely --- if the initial divergence  is ``not too negative'' as in \eqref{eq:CT2D}, and the initial spectral gap is ``not too large'' as in \eqref{eq:etaCT}, then global regularity persists for all time.  In particular, since $\etaS=\sqrt{(\pa_1u_1-\pa_2u_2)^2 + (\pa_1u_2+\pa_2u_1)^2}$ we find that both \eqref{eq:CT2D},\eqref{eq:etaCT} hold if 
\[
|\pa_ju_i (x,0)| \leq \frac{1}{4\sqrt{2}}m_0\phi_\infty.
\]
\end{rmk}

\begin{rmk}[on the finite variation \eqref{V_0}] Observe that \eqref{V_0} places a restriction on the size of $V_0$ even in the case of unconditional flocking, $|\phi|_1=\infty$. Specifically, recall that $V_0$ 
dictates the maximal diameter of the flock in \eqref{eq:Dinf} and thus, \eqref{V_0} amounts to
\begin{equation}\label{eq:finiteD}
\int_{D_0}^{\diam_\infty} \phi(s)ds \leq  \frac{\phi^2(\diam_\infty)}{4\,{\mathop{\max}}_{s\leq \diam_\infty}|\phi'(s)|}.
\end{equation}
Since the term on the left is increasing  while the term on the right is decreasing as functions of  $\diam_\infty$, it follows that \eqref{eq:finiteD} is satisfied for diameters $\diam_\infty$ up to some maximal finite size, that is --- the condition made in \eqref{V_0} is met for finite 
$\displaystyle V_0=m_0\int^{\diam_\infty}\hspace*{-0.2cm}\phi(s)ds$ depending on the influence function $\phi$. This finite restriction on $V_0$ can probably be improved, but unlike the one-dimensional case it cannot be completely removed. In fact, since $V_0 \leq (\mu_2(S_0) + \omega_0)\diam_\infty$, the bound sought in \eqref{V_0} places a purely two-dimensional restriction on the size of \emph{initial vorticity}. 
\end{rmk}

\begin{rmk}[on the finite horizon] Observe that in the case of alignment with a finite horizon, the critical threshold \eqref{eq:CT2D} requires that $\diver{\u_0}(x)\geq 0$ for $\text{dist}\{x, \supp\{\rho_0\}\}>D_\infty$. This is precisely the critical threshold  condition which rules out finite time blow-up in the pressure-less equations \cite{Ta2017}, which is satisfied when prescribing far-field constant velocity \eqref{eq:BC}. In this case,  the critical threshold \eqref{eq:CT2D} needs to be verfied within the finite horizon $\text{dist}\{x, \supp\{\rho_0\}\} < D_\infty$.
\end{rmk}

\begin{proof} 
Our purpose is to show that the derivative $\{\pa_ju_i\}$ are uniformly bounded. We proceed in four steps.

\medskip\noindent
\underline{Step \#1} --- the dynamics of $\diver{\u}+\phi*\rho$. Differentiation of \eqref{The basic equation} implies that  the $2\times 2$ velocity gradient matrix, $M_{ij}:=\pa_j u_i$, satisfies

\bel\label{eq:CSM}
M_t+\u\cdot\nabla M + M^2 
=-(\phi*\rho) M+R, \qquad R_{ij}:=\pa_j\phi *(\rho u_i) -u_i\pa_j \phi*\rho.
\eel
The entries of the residual matrix $\{R_{ij}\}$ can bounded by the commutator estimate \cite[proposition 4.1]{TT2014} in terms of $\displaystyle V(t)=\sup_{\text{supp}(\rho)} \hspace*{-0.05cm}|u_i(x,t)-u_i(y,t)| \leq V_0e^{-\kappa t}$,
\[
|R_{ij}|= \left|\int_{\rr^n}\pa_j \phi(|x-y|)(u_i(y,t)-u_i(x,t))\rho(y,t)\dy\right|\leq |\phi'|_\infty m_0V_0e^{-\kappa t}, \qquad \kappa=m_0\phi_\infty.
\]
The first step is to bound the divergence: taking the trace of \eqref{eq:CSM} we find that   $\d:=\nabla\cdot \u$ satisfies
\[
\d_t+ \u\cdot\nabla \d + \trace{M^2} = - (\phi*\rho)\d +\trace{R}.
\] 
Expressed in terms of the material derivative along particle path, $X':=(\pa_t+\u\cdot\nabla)X$, we have $\d' + \trace{M^2} = - (\phi*\rho)\d +\trace{R}$.
We now make  a key observation that $\trace{R}$ is in fact an exact  derivative along particle path. Indeed, as in \cite{CCTT2016} we invoke the mass equation, 
\[
\trace{R}= \phi* \nabla\cdot (\rho \u) - \u\cdot\nabla \phi*\rho = 
-(\phi* \rho)_t -\u\cdot\nabla \phi*\rho = -(\phi*\rho)',
\]
 and we end up with
\begin{equation}\label{eq:junc}
(\d+\phi*\rho)'+ \trace{M^2} = -(\phi*\rho)\d.
\end{equation}
To proceed, we  express 
$\displaystyle \trace{M^2} \equiv \frac{\d^2+\etaM^2}{2}$ in terms of the \emph{spectral gap},  $\etaM:=\lam_2(M)-\lam_1(M)$, associated with the eigenvalues of $M$,
\bel\label{eq:CSd}
(\d+\phi*\rho)' = -\hf\etaM^2 -\hf\d(\d+2\phi*\rho).
\eel
We need to follow the dynamics of the spectral gap $\eta_M$. To this end, one may try to use the \emph{spectral dynamics} associated with $M$, \cite{LT2002}: by \eqref{eq:CSM} the $\lam_i$'s satisfy
\[
\lam'_i+\lam_i^2 =- (\phi*\rho)\lam_i + \langle \bell_i, R \br_i\rangle, \qquad i=1,2,
\]
where $\{\bell_i,\br_i\}$ are the left and right eigenvectors associated with $\lam_i$, normalized such that $\langle \bell_i,\br_i\rangle=1$.  
Taking the difference of these two equations shows that the spectral gap $\etaM=\lam_2-\lam_1$, satisfies the transport equation
\[
\etaM' + (\d+\phi*\rho)\etaM = \langle \bell_2, R \br_2\rangle - \langle \bell_1, R \br_1\rangle.
\]
Here one faces the difficulty which arises with the term on the right, namely ---  even with the control of the entries $\{R_{ij}\}$, we may still encounter an ill-conditioned $M$ with $|\bell_i|\cdot|\br_i| \gg1$ so that the magnitude of this term is left unchecked. To circumvent this difficulty, we proceed along the lines argued in \cite{Ta2017}: we split $M$ into its symmetric and antisymmetric parts $M=S+\Omega$ and  use the identity\footnote{Equating the trace of $M^2$ with that  of $S^2+\Omega^2 +S\Omega + \Omega S$ we find $\trace{M^2}=\trace{S^2}-2\omega^2$. Using $\trace{X^2}= \hf({\d^2+\eta_{{}_X}^2})$ with $X=M$ on the left and $X=S$ on the right implies \eqref{eq:etasym}.} 
\begin{equation}\label{eq:etasym}
\etaM^2\equiv\etaS^2-4\omega^2, \qquad M=S+\Omega, \ \ \Omega:=\left(\begin{array}{cc}0 & -\omega\\ \omega & 0 \end{array}\right),
\end{equation}
where $\omega$ is the \emph{scaled} vorticity\footnote{The use of such scaling simplifies the computation below.} $\omega=\frac{1}{2}(\pa_1u_2-\pa_2u_1)$. Expressed in terms of $\etaS$, the trace dynamics \eqref{eq:CSd} now reads
\[
(\d+\phi*\rho)' = \hf(4\omega^2 -\etaS^2)-\hf \d(\d+2\phi*\rho).
\]
This calls for the introduction  of the new ``natural'' variable $\e=\d+\phi*\rho$, satisfying
\bel\label{eq:CSe}
 \e' = \hf \left((\phi*\rho)^2 +4\omega^2 -\etaS^2 -\e^2\right). 
 \eel 
 Our purpose is to show that $\{x \ | \ \e(x,t) \geq0\}$  is invariant region of the dynamics \eqref{eq:CSe}.
 
\medskip\noindent  
\underline{Step \#2} --- bounding the spectral gap $\etaS$. Consider the dynamics of the symmetric part of \eqref{eq:CSM}
\[
S' + S^2 =\omega^2{\mathbb I}_{2\times 2}  -(\phi*\rho)S  + R_{\text{sym}}, \qquad R_{\text{sym}}=\hf(R+R^\top), \quad  
\]
The spectral dynamics of its eigenvalues, $\mu_2(S)\geq \mu_1(S)$, is governed by
\begin{equation}\label{eq:mui}
\mu'_i + \mu_i^2 =\omega^2 -(\phi*\rho)\mu_i +\big\langle \bs_i, R_{\text{sym}}\bs_i\big\rangle
\end{equation}
driven by  the \emph{orthonormal} eigenpair $\{\bs_1,\bs_2\}$ of the symmetric $S$.
Taking the difference, we find that $\etaS:=\mu_2(S)-\mu_1(S)\geq 0$ satisfies,
\begin{equation}\label{eq:etaS}
\etaS' + \e\etaS = q, \qquad \e=\d+\phi*\rho.
\end{equation}
 This is the \emph{same} dynamics  found with  $\etaM$ except that the different residual on the right of \eqref{eq:etaS} given by
\[
q:= \big\langle \bs_2, R_{\text{sym}}\bs_2\big\rangle - \big\langle \bs_1, R_{\text{sym}}\bs_1\big\rangle,
\]
is  now controlled by the size of $\{R_{ij}\}$: since  $\bs_i$ are normalized,
\begin{equation}\label{eq:res}
|q(\cdot,t)| \leq 2 \max_{ij} |R_{ij}(\cdot,t)| \leq 2|\phi'|_{\infty} m_0 V_0 e^{-\kappa t}, \qquad \kappa=m_0\phi_\infty.
\end{equation}
Hence, as long as $\e(\cdot,t)$ remains positive then $\etaS$ remain uniformly bounded
\begin{equation}\label{eq:etaSbd}
|\etaS(x,t)| \leq \max_x |\etaS(x,0)| + 2 \frac{|\phi'|_\infty}{\phi_\infty}V_0 
< \max_x |\etaS(x,0)| +  \frac{1}{2} m_0\phi_\infty < m_0\phi_\infty
\end{equation}
The first inequality on the right follows from integration of \eqref{eq:etaS}-\eqref{eq:res}; the second follows from the $V_0$-bound in \eqref{V_0} and the third from the assumed bound on $\etaSz$ in \eqref{eq:etaCT}.

\medskip\noindent
\underline{Step \#3} --- the invariance of $\e(\cdot, t)\geq 0$ . We return to \eqref{eq:CSe}: expressed in terms of $c(x,t):=\sqrt{(\phi*\rho)^2-\etaS^2}$ we have  
\begin{equation}\label{eq:eq}
\e'\geq \hf \left(c^2(x,t) - \e^2\right), \qquad c(x,t)=\sqrt{(\phi*\rho)^2- \etaS^2}.
\end{equation}
Observe that $c(\cdot)$ is well-defined in $\rr$: the upper-bound \eqref{eq:etaSbd} and the lower-bound $\phi*\rho \geq m_0\phi_\infty$ imply that as long as $\e\geq0$, the right term on the right of \eqref{eq:eq} remains boundedly positive 
\[
c(x,t) \geq \sqrt{m^2_0\phi_\infty^2- \max_x \etaS^2(x,t)} \geq \cm >0. 
\]
Since $\e'  \geq \hf(\cm^2-\e^2)= \hf(\cm-\e)(\cm+\e)$,  it follows that $\e$ is increasing whenever $\e\in (-\cm,\cm)$ and in particular, if $\e_0\geq0$ then $\e(x,t)$ remains positive at later times. Thus, if the initial data are \emph{sub-critical} in the sense that \eqref{eq:CT2D} holds  
\[
\e_0= \diver{\u_0}(x)+\phi*\rho_0(x) >0,
\]
then $\e(\cdot,t)\geq 0$ and $\etaS(\cdot,t)$ remains bounded.

\medskip\noindent
\underline{Step \#4} --- an upper-bound of $\e(\cdot,t)$. The lower-bound  $\e \geq 0$ implies that the vorticity is bounded. Indeed, the anti-symmetric part of \eqref{eq:CSM} yields that the vorticity $\omega=\hf\trace{JM}$ satisfies  
\begin{equation}\label{eq:vorticity}
\omega' +\e\omega= \hf\trace{JR}, \qquad J=\left(\begin{array}{cc} 0 & -1 \\ 1 &0\end{array}\right)
\end{equation}
hence 
\begin{equation}\label{eq:vort}
|\omega|'\leq  -\e|\omega| + \hf|q|, \qquad |q(\cdot,t)|\leq 2|\phi'|_\infty m_0V_0e^{-\kappa t}, \quad \kappa=m_0\phi_\infty,
\end{equation}
and we end up with same upper-bound on $\omega$ as with $\etaS$,
\begin{equation}\label{eq:vorbd}
|\omega(x,t)| \leq \omega_{\text{max}}, \qquad \omega_{\text{max}}:= \max_x |\omega_0| +  \hf m_0\phi_\infty.
\end{equation}
Returning to \eqref{eq:CSe} we have (recall $\phi\leq 1$)
\[
\e' \leq \hf \Big((\phi*\rho)^2 +4\omega^2 -\e^2\Big) \leq \hf \Big(m_0^2+ 4\omega^2_{\text{max}} -\e^2\Big),
\]
which implies that $\e(x,t) \leq \e_{\text{max}} <\infty$.
The uniform bound on $\e$ implies that $\diver{u}$ is uniformly bounded,
$|\diver{\u}| \leq |\e|_\infty+|\phi*\rho|_\infty \leq \e_{\text{max}}+ m_0$, and together with the bound on the spectral gap \eqref{eq:etaSbd}, it follows that the symmetric part $\{S_{ij}\}$ is bounded. 
Finally, together with the vorticity bound \eqref{eq:vorbd} it follows that $\{\pa_j u_i\}$ are uniformly bounded which completes the proof. 
\end{proof}

\begin{rmk}
Observe that the region of sub-critical configuration leading global regularity becomes \emph{larger} for $|\omega_0| \gg1 $ in agreement with the statements made in \cite{LT2004, CT2008} that rotation prevents or at least delays finite-time blow-up. Specifically, if $|\omega_0(\cdot)| \geq \omega_{\text{min}}>0$ then one can set a larger lower barrier $c= \sqrt{(\phi*\rho)^2 +4\omega^2_{\text{min}}-\etaS^2}$ in \eqref{eq:eq} leading to the improved threshold
$\diver{\u_0} > -\phi*\rho_0-\omega_{\text{min}}$. In  particular, if $\omega$ is large enough so that $4\omega^2-\eta^2_S>0$, that is --- if $M$  has complex-valued eigenvalues, then the invariance of the positivity of $\e$ follows at once from  the fact that  \eqref{eq:CSe}   is dominated equation  by  $\e' \geq \hf \left((\phi*\rho)^2  -\e^2\right)$.  As in the  2D restricted Euler-Poisson equations \cite{LT2003}, the difficulty lies with the case of real eigenvalues. 
\end{rmk}

\begin{rmk}
The proof of theorem \ref{thm:CS} reveals two main aspects. First, the \emph{commutator structure} of the alignment term on the right of \eqref{CS}${}_2$,  expressed as
$[\phi*, u](\rho)$, leads to the `residual terms' $R_{ij}$ with exponentially decaying bound. The role of commutator structure was highlighted in our recent work \cite{ST2016}. Second, the use of spectral dynamics, \cite{LT2002,LT2003,LL2013}, to trace the propagation of regularity for the remaining, non-residual terms
in \eqref{eq:CSM}.   
\end{rmk}

\subsection{Fast alignment}\label{sec:CS-flock}
We extend the one-dimensional arguments of \cite{ST2016} which show that exponentially rapid convergence  towards a \emph{flocking state}, consisting of  a constant 2-vector velocity  $\{\bar{\u}\in \mathbb{R}^2$ and a traveling density profile $\bar{\rho}(x,t)=\rho_\infty(x-t\bar{\u})\}$.  We  only indicate the main aspects in the passage  to the present system.  We start by noting that the positivity of $\e$ implies more than the mere boundedness of the spectral gap $\eta_S$ and the vorticity $\omega$. Indeed, \eqref{eq:etaS} and \eqref{eq:vort} imply that these quantities  follow the  exponential decay of $q$ in \eqref{eq:res}
\[
|\eta_S(\cdot,t)|_\infty + |\omega(\cdot,t)|_\infty \lesssim e^{-\kappa t}.
\]
This shows that modulo rapidly decaying error terms $E(t)$ of order $E(t) \lesssim e^{-\kappa t}$,  equation \eqref{eq:CSe} which governs $\e$ takes the form
\[
\e_t +\u\cdot\nabla\e=\frac{1}{2}\left(\h^2-\e^2\right) + E(t), \qquad \h:=\phi*\rho
\] 
Moreover, convolving the mass equation with $\phi$ we find
\begin{equation}\label{eq:convrho}
\h_t +\u\cdot\nabla \h = \int \nabla\phi(|x-y|)\cdot (\u(x,t)-\u(y,t))\rho(y,t)\dy.
\end{equation}
Observe that the quantity on the right of rapidly decaying, being upper-bounded by  $ \lesssim |\phi'|_\infty V(t) \lesssim e^{-\kappa t}$. Hence, the difference $\d=\e-\h$ satisfies
\[
\d_t +\u\cdot\nabla\d = -\frac{1}{2}(\e+\h)\d + E(t).
\] 
The positivity of $\e+\h$ then implies the rapid decay of the divergence, $|\diver{\u}(\cdot,t)|_\infty \lesssim e^{-\kappa t}$. The exponential decay of the divergence, the vorticity and the spectral gap imply that $|\partial_ju_i(\cdot,t)|_\infty \lesssim e^{-\kappa t}$. Let $\bar{\u}$ be a large-time limiting value of $\u(\cdot,t)$. The mass equation reads
\[
\rho_t + \bar{\u}\cdot\nabla\rho= -\d\rho +(\bar{\u}-\u)\cdot\nabla\rho.
\]
The term on the right is rapidly decaying because $\d$ and $(\bar{\u}-\u)$ are, and one concludes along the lines of \cite{ST2017}, that there exists a traveling density profile such that $\rho(x,t)-\rho_\infty(x-t\bar{\u}) \rightarrow 0$.

\section{Motsch-Tadmor  hydrodynamics: global regularity and fast alignment}
In this section, we study the flocking hydrodynamics which arises from MT model \eqref{eq:fhd} with $\kappa= \phi_\infty$. We begin by recalling the one-dimensional case
\bel\label{MT1D}\ba
&\rho_t+(\rho u)_x=0, \qquad (x,t)\in (\rr,\rr_+)\\
&u_t+uu_x=\int \frac{\phi(|x-y|)}{(\phi*\rho)(x,t)}(u(y,t)-u(x,t))\rho(y,t)\dy.
\ea\eel
System \eqref{MT1D}  was recently studied in \cite{BRSW2015}, as the hydrodynamic description for agent-based model of ``emotional contagion'', and in \cite{GG2017} in the context of stable swarming.
In \cite{CCTT2016} it was proved that  (\ref{MT1D}) has a global classical solution for sub-critical  initial data such that 
\begin{equation}\label{eq:MT1D}
\pa_x u_0(x)\geq -\sigma_+(V_0) \ \mbox{ for all} \ \ x\in \rr,
\end{equation}
for a certain critical curve $\sigma_+ \geq 0$.
We now make a precise statement of the critical threshold  for both the one - and two-dimensional MT model.

\begin{thm}
[Critical threshold for 2D Motsch-Tadmor hydrodynamics]\label{thm:MT}
\noindent
Consider the two-dimensional MT model in $(x,t)\in (\rr^2,\rr_+)$,
\bel\label{MT2D}
\left\{
\ba
&\rho_t+\na\cdot(\rho \u)=0,\\
&\u_t+\u\cdot\na \u=\int a(x,y,t)(\u(y,t)-\u(x,t))\rho(y,t)\dy, \qquad  a(x,y,t):=\frac{\phi(|x-y|)}{(\phi*\rho)(x,t)},
\ea
\right.
\eel
subject to initial conditions $(\rho_0, \u_0)\in (L^1, W^{1,\infty}(\rr^2))$,
with compactly supported density, $D_0 <\infty$ and initial velocity of finite variation
\bel\label{V_1}
V_0\leq m_0\cdot\min\left\{|\phi|_{1},\frac{\phi^2_\infty}{4|\phi'|_\infty(1+2\phi_\infty)}\right\}, \qquad \phi_\infty=\phi(\diam_\infty).
\eel

\noindent
Assume that the following critical threshold condition holds.\newline
 (i) The initial velocity divergence satisfies
\bel\label{eq:MT2D}
\diver{\u_0}(x) \geq -1 \quad \mbox{for all} \quad x\in \rr^2.
\eel
(ii) Then the initial \emph{spectral gap} $\etaSz:=\mu_2(S_0)-\mu_1(S_0)$ is bounded
\begin{equation}\label{eq:etaMT}
\max_x \big|\etaSz(x)\big| \leq \frac{1}{2}, \qquad \etaS=\mu_2(S(x,t))-\mu_1(S(x,t)).
\end{equation}  
Then the class of such sub-critical initial conditions \eqref{eq:MT2D},\eqref{eq:etaMT} give rise to a classical solution $(\rho(t),\u(t) \in C(\rr^+;L^\infty(\rr^2 ))\times C(\rr^+;\dot{W}^{1,\infty}(\rr^2))$ with large time hydrodynamics flocking behavior \eqref{S, V exp-control} $\displaystyle \mathop{\max}_{x\in \text{supp}(\rho)} |\u(x,t)-\u(y,t)| \lesssim e^{-\kappa t}$.
 \end{thm}
\begin{rmk}
In the case of  finite horizon alignment encoded in \eqref{eqs:horizon} with $\al=\phi*\rho$,  the critical thresholds \eqref{eq:MT2D},\eqref{eq:etaMT} can be restricted to the finite set $\text{dist}\{x,\supp\{\rho_0\}\}$.
\end{rmk}
\begin{proof}
As before, we trace the dynamics of $M=\pa_ju_i$,
\bel\label{eq:MTM}
M_t+\u\cdot\nabla M + M^2  =- M+R, 
\eel
where the entries of the residual matrix $\{R_{ij}\}$ are given by
\[
R_{ij}(x,t):=\int_{y\in \rr^2}\pa_j a(x,y,t)(u_i(y,t) -u_i(x,t))\rho(y,t)\dy, \qquad a(x,y,t)=\frac{\phi(|x-y|)}{(\phi*\rho)(x,t)}
\]
Expressed in terms of the operator $A(w):=\int_y a(x,y,t)w(y)dy$, the entries of $R$ have the commutator structure $R_{ij}= \pa_j [A, u_i](\rho) $ which can be estimated by the commutator bound \cite[proposition 7.1]{TT2014} in terms of $V(t)=\sup_{\text{supp}(\rho)} \hspace*{-0.05cm}|u_i(x,t)-u_i(y,t)|$,
\[
|R_{ij}(x,t)| = \big|\pa_j [A, u_i](\rho) \big|\leq  \frac{|\phi'|_\infty}{\phi_\infty}V_0e^{-\kappa t}, \qquad \kappa=\phi_\infty.
\]
We now proceed as before. As a first step, we follow the dynamics of the $\d=\diver{\u}$: taking the trace of ${(\ref{eq:MTM})}$ we find 
\begin{equation}\label{eq:rb}
\d' + \hf (\d^2+\etaS^2) = \omega^2-\d + r, \qquad r:=\trace{R} \leq 2\frac{|\phi'|_\infty}{\phi_\infty}V_0.
\end{equation}
This calls for the introduction of a new variable $\e:=\d+1$ where the last equation recast into the Riaccti's form
\begin{equation}\label{eq:MTe}
\e'= \hf \Big(1-\etaS^2 +2r -\e^2\Big) +\omega^2. 
\end{equation}
Our purpose si to show that the $\{x \ | \ \e(x,t)\geq0 \}$ is invariant of the dynamics \eqref{eq:MTe} and to this end we need to bound the spectral gap $\etaS$.

The second step is to follow the spectral dynamics associated with the symmetric part of \eqref{eq:MTM}
\[
\mu'_i(S) + \mu^2_i(S) = \omega^2 -\mu_i(S) +\big\langle \bs_i,R_{\text{sym}}\bs_i\big\rangle.
\]
Taking the difference and recalling that $\bs_i$ are the normalized eigenvectors of $S$ we find the dynamics of the spectral gap,
\begin{equation}
\etaS' + \e \etaS = q, \qquad |q| \leq 2\max|R_{ij}(x,t)| \leq 2\frac{|\phi'|_\infty}{\phi_\infty}V_0 e^{-\kappa t}.
\end{equation}
It follows that as long as $\e(\cdot,t)$ is positive then
\begin{equation}\label{eq:MTetab}
|\etaS(x,t)| \leq \max_x |\etaSz(x)| + 2\frac{|\phi'|_\infty}{\phi^2_\infty}V_0 <\hf,
\end{equation}
and therefore $c:=\sqrt{1-\etaS^2+2r}$ has the lower bound $c(x,t)\geq c_{\text{min}}>0$, where
\[
\max_x |\etaSz(x)| + \Big(2\frac{|\phi'|_\infty}{\phi^2_\infty} + 
 4\frac{|\phi'|_\infty}{\phi_\infty}\Big)V_0 \leq 1-c^2_{\text{min}} <1
\]
This inequality follows from the assumed bounds on $V_0$ in \eqref{V_1} and on the initial spectral gap \eqref{eq:etaMT}, and the bound of $r$ in \eqref{eq:rb}. As a final step, we  return to \eqref{eq:MTe} to find, $\e' \geq \hf(c^2_{\text{min}}-\e^2)$,
which guarantees that if the critical threshold \eqref{eq:MT2D} holds, i.e., if $\e_0\geq0$ then $\e(x,t)\geq0$ at later time. Moreover, since $\e(\cdot, t)\geq0$, the vorticity equation,
 $\omega'+\e\omega = \hf\trace{JR}$,
shows that $|\omega(\cdot,t)|$ remains bounded in terms of $\max_x|R_{ij}(x,t)| \lesssim r_{\text{max}} <\infty$. The transport equation \eqref{eq:MTe} implies 
\[
\e' \leq \hf\Big(1+2r  + 2\omega^2 - \e^2\Big) \leq \hf\Big(\frac{3}{2} + 2 \omega^2_{\text{max}}-\e^2\Big),
\]
 and a uniform  upper-bound of $\e(\cdot, t) \leq \e_{\text{max}} <\infty$ follows.
\end{proof}

\begin{rmk}
In the one-dimensional case, $\etaS=\omega\equiv 0$ and the dynamics of $\e=\d+1$ in \eqref{eq:MTe} simplifies into $\e'=\hf(1+2r -\e^2)$. Hence, the variation bound \eqref{V_1} can be related to 
\[
V_0< m_0\min\left\{|\phi|_{1}\,, \frac{1}{4} \frac{\phi_\infty}{|\phi'|_\infty}\right\}
\]
so that $1+2r \geq c_{\text{min}}>0$ and  $\e' > \hf(c_{\text{min}}-\e^2)$ implies global smoothness under the critical threshold condition $\pa_xu_0(x) \geq -1$.
\end{rmk}

\begin{rmk}
One can follow the argument in section \ref{sec:CS-flock} to conclude that the same rapid alignment holds for MT model. Indeed, the MT model enhances the convergence \underline{rate} towards a limiting flocking state. 
\end{rmk}

\end{document}